\documentclass[final]{siamltex}
\usepackage{amsmath}
\usepackage{amssymb}
\usepackage{xspace}

\newenvironment{Theorem}{\begin{theorem}}{\end{theorem}}
\newenvironment{Proposition}{\begin{proposition}}{\end{proposition}}
\newenvironment{Lemma}{\begin{lemma}}{\end{lemma}}
\newenvironment{Corollary}{\begin{corollary}}{\end{corollary}}
\newenvironment{Definition}{\begin{definition}}{\end{definition}}
\newenvironment{Proof}{\begin{proof}}{\end{proof}}
\newtheorem{Example}{Example}[section]
\newtheorem{Remark}{Remark}[section]
\newtheorem{Conjecture}{Conjecture}[section]

\newcommand{\two}{\mathbb{F}_2}
\newcommand{\sub}[1]{[ #1 ]}
\newcommand{\xor}{\oplus}
\newcommand{\sdif}{\mathop{\mathrm{\Delta}}}
\newcommand{\vertexrem}{\setminus}
\newcommand{\setminor}{\setminus}
\newcommand{\dual}{\mathrel{\bar*}}
\newcommand{\minc}{\mathrm{minc}}
\newcommand{\maxc}{\mathrm{maxc}}
\newcommand{\dmatroid}{delta-matroid\xspace}
\newcommand{\dmatroids}{delta-matroids\xspace}
\newcommand{\vfsafe}{vf-safe\xspace}

\renewcommand{\emptyset}{\varnothing}

\title{Nullity and Loop Complementation for Delta-Matroids}
\author{Robert Brijder\thanks{Hasselt University and Transnational University of Limburg, Belgium, \texttt{robert.brijder@uhasselt.be}} \and Hendrik Jan Hoogeboom\thanks{Leiden Institute of Advanced Computer Science,
Leiden University, The Netherlands, \texttt{hoogeboom@liacs.nl}}}

\begin{document}

\maketitle

\begin{abstract}
We show that the symmetric difference distance measure for set
systems, and more specifically for delta-matroids, corresponds
to the notion of nullity for symmetric and skew-symmetric
matrices. In particular, as graphs (i.e., symmetric matrices
over GF(2)) may be seen as a special class of delta-matroids,
this distance measure \emph{generalizes} the notion of nullity
in this case. We characterize delta-matroids in terms of
equicardinality of minimal sets with respect to inclusion (in addition
we obtain similar characterizations for matroids). In this way,
we find that, e.g., the delta-matroids obtained after loop
complementation and after pivot on a single element together
with the original delta-matroid fulfill the property that two
of them have equal ``null space'' while the third has a larger
dimension.
\end{abstract}



\section{Introduction}
A \emph{set system} is a tuple $M = (V,D)$ with $V$ a finite
set, called the \emph{ground set}, and $D$ a family of subsets
of $V$. Set system $M$ is called \emph{proper} if $D \not=
\emptyset$. Let $X \subseteq V$.
The \emph{pivot}
(also called \emph{twist}) of $M$ on $X$, denoted by $M * X$,
as $(V,D * X)$, where $D * X = \{Y \sdif X \mid Y \in D\}$
\cite{bouchet1987}  (here $\sdif$ denotes symmetric
difference).
We denote by $\min(M)$ (and $\max(M)$,
resp.) the set system with ground set $V$ consisting of the
minimal (maximal, resp.) sets with respect to set inclusion of a set
system $M$. A \dmatroid is a proper set system $M$ that
satisfies the \emph{symmetric exchange axiom}:
For all $X,Y \in D$ and all
$u \in X \sdif Y$,
there is a $v \in X \sdif Y$ (possibly $v=u$) such that $X
\sdif \{u,v\} \in D$ \cite{bouchet1987}.

The main results of this paper are described below.
We characterize first the notion of a \dmatroid.

%
{\scshape Theorem~\ref{thm:isodistant_dmatroid_eq}.}
\begin{em}
Let $M$ be a proper set system. Then $M$ is a \dmatroid iff for
each $X \subseteq V$, the sets in $\min(M*X)$ have equal
cardinality.
\end{em}

We will almost exclusively work with this
characterization rather than directly using the symmetric
exchange axiom. Related to Theorem~\ref{thm:isodistant_dmatroid_eq},
we obtain novel characterizations of a matroid described by
its bases and its independent sets. Below is the
characterization of a matroid described by its independent
sets.

{\scshape Corollary~\ref{cor:char_matroids_indep_sets}.}
\begin{em}
Let $M$ be a proper set system. Then $M$ is a matroid described
by its independent sets iff both (1) for each $Y \in M$ and $Y'
\subseteq Y$, we have $Y' \in M$, and (2) for each $X \subseteq
V$, the sets in $\min(M*X)$ have equal cardinality.
\end{em}

Given a set system $M$ and a set $X$, the distance $d_M(X) =
\min(\{|X \sdif Y| \mid Y \in M\})$ is the minimal cardinality
of the symmetric difference of $X$ and the sets in $M$. It
turns out the distance behaves well under deletion of elements
from the ground set provided we consider \dmatroids, cf.
Theorem~\ref{thm:iso_invariant_minors}.

For a set system $M = (V,D)$ and $v \in V$, we define
\emph{pseudo-deletion} of $M$ on $v$, denoted by $M
\widehat\vertexrem v$, as $(V,D')$, where $D' = \{X \cup \{v\}
\mid X \in D, v \not\in X\}$. Moreover, we define \emph{loop
complementation} of $M$ on $v$, denoted by $M + v$, as
$(V,D'')$, where $D'' = D \sdif D' = D \sdif \{X \cup \{v\}
\mid X \in D, v \not\in X\}$ \cite{BH/PivotLoopCompl/09}. Loop
complementation is well motivated: it generalizes the loop
complementation for graphs (this is recalled in
Section~\ref{sec:null_graphs}). We derive the following
property of \dmatroids in relation to pivot and loop
complementation.

{\scshape Theorem~\ref{thm:nullspace_plus_op_setsystem_max}.}
\begin{em}
Let $M$ be a \dmatroid, and $v \in V$ such that $M+v$ is a
\dmatroid. Then $\max(M)$, $\max(M*v)$, and $\max(M+v)$ are
matroids (described by their bases) such that precisely two of
the three are equal, to say $M_1$. Moreover, the rank of the
third $M_2$ is one smaller than the rank of $M_1$ and $M_2
\widehat\vertexrem v = M_1$.
\end{em}

We also formulate a ``$\min$ counterpart'' of
Theorem~\ref{thm:nullspace_plus_op_setsystem_max}, cf.\
Theorem~\ref{thm:null_plus_op_setsystem}, which involves the
dual pivot operation instead of the loop complementation
operation.

Given a $V \times V$ matrix $A$ (the columns and rows of $A$
are indexed by finite set $V$), we denote by $A\sub{X}$ the
principal submatrix of $A$ induced by $X$ and we define the set
system $\mathcal{M}_A = (V,D_A)$ with $D_A = \{ X \subseteq V
\mid A\sub{X} \mbox{ is nonsingular}\}$. If $A$ is
skew-symmetric (i.e., $A^T = -A$ where $A^T$ denotes the
transpose of $A$) or symmetric, then $\mathcal{M}_A$ is a
\dmatroid \cite{bouchet1987}. We show that nullity of a
principal submatrix $A[X]$ corresponds to distance in the
associated \dmatroid $\mathcal{M}_A$. This is the main
motivation for considering distance, as it allows us to carry
over results of distances in \dmatroids to nullity values of
symmetric or skew-symmetric matrices and graphs in particular.

{\scshape Theorem~\ref{thm:respr_setsystem_isodistant_skew_symm}.}
\begin{em}
Let $A$ be a $V \times V$ symmetric or skew-symmetric matrix
(over some field). Then $d_{\mathcal{M}_A}(X) = n(A\sub{X})$
for each $X \subseteq V$.
\end{em}

It is known that $\mathcal{M}_A*X$, for any $V \times V$ matrix
$A$ and any set $X$ in $\mathcal{M}_A$, is equal to
$\mathcal{M}_{A*X}$ where $A*X$ is the principal pivot
transform of $X$ on $A$ (see
\cite{tucker1960,Tsatsomeros2000151} for the definition of this
notion). Hence there is a close connection between the linear
algebra of principal pivot transform and the combinatorics of
pivot on set systems.

The \dmatroid $\mathcal{M}_A$ for a symmetric or skew-symmetric
matrix $A$ is especially interesting over the binary field
$\two$ (note that the notions of skew-symmetric and symmetric
coincide over $\two$), i.e., in the case where $A$ is a graph
(where loops are allowed). In this case $\mathcal{M}_A$ retains
all information of $A$
--- hence $\mathcal{M}_A$ is a representation of the graph $A$.
It turns out that for a graph the null space (kernel) of its
adjacency matrix is determined by the set of maximal elements
in the associated \dmatroid. In this way we obtain the
following result (we associate a graph $G$ by its adjacency
matrix $A(G)$). For a graph and vertex $v$, $G+v$ denotes loop
complementation (the existence of a loop on $v$ is
complemented), and if $v$ is a looped vertex, then $G*v$
denotes principal pivot transform on $v$, which over $\two$ is
local complementation (the subgraph of the neighbourhood of $v$
is complemented). As usual, we identify vectors indexed by $V$
over $\two$ by subsets of $V$.

{\scshape Theorem~\ref{thm:nullspace_plus_op_graph}.}
\begin{em}
Let $G$ be a graph having a looped vertex $v$. Then $\ker(G)$,
$\ker(G*v)$, and $\ker(G+v)$ are such that precisely two of the
three are equal, to say $K_1$, and the third, $K_2$, is such
that $\dim(K_2) = \dim(K_1) + 1$ and $K_1 = \{ X \in K_2 \mid v
\not\in X\}$.
\end{em}

This result is related to \cite[Lemma~23]{LT/LAA/2011} (there a
graph different from $G*v$ is considered) and
\cite[Theorem~(9.4)]{Bouchet/87/ejc/isotropicsys}, and it can
be seen as an extension of
\cite[Lemma~2]{DBLP:journals/ejc/BalisterBCP02}. In case $G$ is
a circle graph, Theorem~\ref{thm:nullspace_plus_op_graph} is
applicable to the theory of closed walks in $4$-regular graphs,
see e.g. \cite{LT/BinNullity/09}.

It is known from \cite{BH/PivotLoopCompl/09} that the family of
\dmatroids is not closed under loop complementation. We show in
Section~\ref{sec:strong_DM} that the family of binary
\dmatroids is closed under pivot, loop complementation and
deletion of elements from the ground set.

The results given in this paper are crucial in a subsequent
paper on interlace polynomials of \dmatroids
\cite{BH/InterlacePolyDM/10}.

\section{Pivot and Loop Complementation on Set Systems}
First we fix basic notions and terminology. We denote the field
consisting of two elements by $\two$. In this field addition
and multiplication are equal to the logical exclusive-or and
logical conjunction, which are denoted by $\xor$ and $\land$
respectively. By carrying over $\xor$ to sets, we obtain the
symmetric difference operator $\sdif$. Hence for sets $A, B
\subseteq V$ and $x \in V$, $x \in A \sdif B$ iff $(x \in A)
\xor (x \in B)$.

A \emph{set system} (over $V$) is a tuple $M = (V,D)$ with $V$
a finite set, called the \emph{ground set}, and $D$ a family of
subsets of $V$. Let $X \subseteq V$. We define $M \sub{X} =
(X,D')$ where $D' = \{Y \in D \mid Y \subseteq X\}$, and define
$M \setminor X = M \sub{V \setminus X}$. Set system $M$ is
called \emph{proper} if $D \not= \emptyset$. Let $\min(D)$
($\max(D)$, resp.) be the family of minimal (maximal, resp.)
sets in $D$ with respect to set inclusion, and let $\min(M) =
(V,\min(D))$ ($\max(M) = (V,\max(D))$, resp.) be the
corresponding set systems. Also, we denote the family of
minimal sets with respect to \emph{cardinality} by $\minc(D)$, i.e., $X
\in \minc(D)$ iff $X \in D$ and $|X| \leq |Y|$ for all $Y \in
D$. We let $\minc(M) = (V,\minc(D))$ be the corresponding set
system. Similarly, we define $\maxc(M) = (V,\maxc(D))$. Note:
we will also use $\min(E)$ and $\max(E)$ for a finite set $E$
of integers, to denote the smallest and largest, resp., integer
in $E$. We simply write $Y \in M$ to denote $Y \in D$, and for
set system $M' = (V,D')$, $M \subseteq M'$ if $D \subseteq D'$.
We also often write $V$ to denote the ground set of the set
system under consideration. A set system $M$ is called
\emph{equicardinal} if for all $X_1, X_2 \in M$, $|X_1| =
|X_2|$.

Let $M = (V,D)$ be a set system. We define, for $X \subseteq
V$, \emph{pivot} of $M$ on $X$, denoted by $M * X$, as $(V,D *
X)$, where $D * X = \{Y \sdif X \mid Y \in D\}$. The pivot
operation (often called \emph{twist} in the literature) is
often denoted by $M \sdif X$ instead of $M * X$ (see, e.g.,
\cite{bouchet1987}). However, as $D * X$ is of course in
general different from $D \sdif X$, to avoid confusion, we use
$*$ for pivot. We define, for $X \subseteq V$, \emph{loop
complementation} of $M$ on $X$ (the motivation for this name is
from graphs, see Section~\ref{sec:null_graphs}), denoted by $M
+ X$, as $(V,D')$, where $Y \in D'$ iff $|\{ Z \in M \mid Y
\setminus X \subseteq Z \subseteq Y \}|$ is odd
\cite{BH/PivotLoopCompl/09}. In particular, if $X = \{v\}$ is a
singleton, then $D' = D \sdif \{Z \cup \{v\} \mid Z \in D, v
\not\in Z\}$.

For notational convenience we often omit the ``braces'' for
singletons $\{v\}$, and write, e.g., $M+v$, $M*v$, and $M
\vertexrem v$.
Loop complementation and pivot
belong to a class of operations called \emph{vertex flips}, cf. \cite{BH/PivotLoopCompl/09}.
Deletion $M\vertexrem u$ is
also a vertex flip operation (modulo a, for this purpose
irrelevant, difference in ground set). To simplify notation, we
assume left associativity of the vertex flips, and write, e.g.,
$M*u+v$ to denote $(M*u)+v$. Vertex flips turn out commute on
different elements. Therefore, if $u,v \in V$ and $u \not= v$,
then, e.g., $M+u\vertexrem v = M\vertexrem v+u$, $M*u\vertexrem
v = M\vertexrem v*u$, $M*u+v= M+v*u$, and $M+u+v= M+v+u$.
Moreover, it is easy to verify that $M+u \vertexrem u = M
\vertexrem u$.

It has been shown in \cite{BH/PivotLoopCompl/09} that pivot
${}*u$ and loop complementation ${}+u$ on a common element $u
\in V$ are involutions (i.e., of order $2$) that generate a
group isomorphic to $S_3$, the group of permutations on
$3$ elements. In particular ${}+u*u+u = {}*u+u*u$
is the third involution,
called the \emph{dual pivot}, and is denoted by $\dual $. We
have, e.g., ${}+u*u = {}\dual  u+u = {}*u \dual  u$ and ${}*u+u
= {}+u \dual  u = {}\dual   u *u$ for $u \in V$ (these are the
two operations of order $3$). The six operations (including the
identity operation) are called \emph{invertible vertex flips}.

It turns out that, for $X \subseteq V$, $M \dual  X = M +X *X
+X$ is equal to $(V,D')$, where $Y \in D'$ iff $|\{ Z \in M
\mid Y \subseteq Z \subseteq Y \cup X \}|$ is odd. In
particular, if $X = \{v\}$ is a singleton, then $D' = D \sdif
\{Z \setminus \{v\} \mid Z \in D, v \in Z\}$. Equivalently, for
$Y \subseteq V$, if $v \in Y$, then $Y \in M \dual  v$ iff $Y
\in M$, and if $v \not\in Y$, then $Y \in M \dual  v$ iff $(Y
\in M) \xor (Y\sdif\{v\} \in M)$.

Finally, it is observed in \cite{BH/PivotLoopCompl/09} that
$\min(M) = \min(M+X)$. Since $\min(M) = \max(M*V)*V$, we have
similarly $\max(M) = \max(M\dual X)$.

We will often use the results of this section without explicit
mention.

\section{Distance in Set Systems}\label{sect:distance-set-systems}

Let $M$ be a proper set system. For $X \subseteq V$, we define
$d_M(X) = \min(\{|X \sdif Y| \mid Y \in M\})$. Hence, $d_M(X)$
is the smallest distance between $X$ and the sets in $M$, where
the distance between two sets is measured as the number of
elements in the symmetric difference. We will study some
properties of this natural notion, and in particular we
investigate the relation between the values $d_{M\rho}(X)$ for
different invertible vertex flips $\rho$ on a fixed element
$v$.

We set $d_M = d_M(\emptyset)$, the cardinality of a smallest
set in $M$.

\begin{Lemma}\label{lem:dM_nullity_prop}
Let $M$ be a proper set system. Then $d_{M*Z}(X) = d_{M}(X
\sdif Z)$ for all $X,Z \subseteq V$.
\end{Lemma}
\begin{Proof}
$d_{M*Z}(X) = \min(\{|X \sdif Y| \mid Y \in M*Z\}) = \min(\{|X
\sdif (Y\sdif Z)| \mid Y \in M\}) = \min(\{|(X \sdif Z) \sdif
Y| \mid Y \in M\}) = d_{M}(X\sdif Z)$.
\end{Proof}

This basic fact is mainly used to reduce (without loss of
generality) results concerning distance from $X \subseteq V$ in
set systems to distance from the empty set, i.e., the
cardinality of the smallest set in $M$: $d_M(X) =
d_{M*X}(\emptyset) = d_{M*X}$.

As  $\min(M) = \min(M+v)$ we infer that the six different
invertible vertex flips on $v$ result in at most three
different values: $d_{M} = d_{M+v}$, $d_{M*v} = d_{M*v+v}$, and
$d_{M\dual v} = d_{M+v*v}$. By Lemma~\ref{lem:dM_nullity_prop}
this can be extended to distance between an arbitrary $X
\subseteq V$ (instead of $\emptyset$) and $M$. When $v\notin X$
then the three equalities above hold essentially unchanged for
distance from $X$ since vertex flip $\rho$ on $v$ and pivot on
$X$ commute: $d_{M\rho}(X) = d_{M\rho*X}(\emptyset) =
d_{M*X\rho}$. However, when $v\in X$ this commutation no longer
holds, and we have to reconsider the equalities. Writing $X' =
X\sdif \{v\}$ and $M'=M*X'$ we have then
\begin{itemize}
\item
$d_{M}(X) = d_{M'*v} = d_{M'*v+v} = d_{M\dual v}(X)$
\item
$d_{M*v}(X) = d_{M'} = d_{M'+v} = d_{M+v*v}(X)$
\item
$d_{M+v}(X) = d_{M'+v*v} = d_{M'\dual v} = d_{M*v+v}(X)$
\end{itemize}

One easily argues that applying an invertible vertex flip changes
$d_M$ by at most one.

\begin{Lemma} \label{lem:dM_props}
Let $M$ be a proper set system. If $\rho$ is an invertible
vertex flip of $M$ on $v \in V$, then (1) $|d_M - d_{M\rho}|
\leq 1$ and (2) $d_{M\rho} \in \{m,m+1\}$ with $m = \min(\{d_M,
d_{M*v}\})$.
\end{Lemma}
\begin{Proof}
Proof of (1). By the above, we need only to verify the cases
$\rho = {}
* v$ and $\rho ={} \dual  v$. By the definitions of pivot and
dual pivot, for any pair of sets $Z, Z \sdif \{v\} \subseteq
V$, at least one of this pair is in $M$ iff at least one of
this pair is in $M\rho$. Hence the smallest cardinality of a
set in $M$ cannot differ by more than one from the smallest
cardinality of a set in $M\rho$.

Proof of (2). By (1), the result is valid for $d_M$ and
$d_{M*v}$, and it suffices to show that $d_{M\dual v} \ge m$.
The argument we use works for any invertible vertex flip
$\rho$. Let $Z \in \minc(M\rho)$, i.e., $Z \in M\rho$ such that
$|Z| = d_{M\rho}$. We have $Z \in M$ or $Z \sdif \{v\} \in M$
(or both). If $Z \in M$, then $d_{M\rho} =|Z| \geq d_M \geq m$.
If $Z \sdif \{v\} \in M$, then $Z \in M*\{v\}$ and $d_{M\rho}
=|Z| \geq d_{M*v} \geq m$. Hence in both cases we have
$d_{M\rho} \geq m$.
\end{Proof}

By Lemma~\ref{lem:dM_props} the three values $d_M$, $d_{M*v}$,
and $d_{M\dual v}$ cannot be all different. As
$d_M(X\sdif\{v\}) = d_{M*\{v\}}(X)$ we also have, for $v \in
V$, $|d_M(X) - d_M(X\sdif\{v\})| \leq 1$.

We obtain now a result for $M\dual v$ assuming $d_M \not= d_{M*v}$.
\begin{Theorem} \label{thm:eq_m_dualp_setsystem}
Let $M$ be a proper set system, and $v \in V$.
We have
$$
\minc(M\dual v) = \begin{cases} \minc(M) & d_M < d_{M*v} \\
\minc(M*v) & d_{M*v} < d_M
\end{cases}.$$
In either case, the elements of $\minc(M\dual v)$ do not
contain $v$. In particular, if $d_M \not= d_{M*v}$, then
$d_{M\dual v} = m$ with $m = \min(\{d_M, d_{M*v}\})$.
\end{Theorem}
\begin{Proof}
We may assume without loss of generality that $d_M < d_{M*v}$.
Indeed, if $d_M > d_{M*v}$, then consider $M' = M*v$. We have
$\minc(M \dual v) = \minc(M' * v \dual v) = \minc(M'\dual v +
v) = \minc(M'\dual v)$ and $d_{M'} < d_{M'*v}$.

Assume therefore that $d_M < d_{M*v}$. Let $Z \in \minc(M)$.
Note that $v\notin Z$ as otherwise $Z \sdif \{v\}$ would be a
smaller set in $M*v$. Moreover $Z \notin M*v$, as otherwise
$d_{M*v} \leq d_{M}$. By the definition of dual pivot, $v
\notin Z$ implies that $Z\in M\dual v$ iff exactly one of $Z$
and $Z\cup\{v\}$ is in $M$. Hence $Z \in M\dual v$. By
Lemma~\ref{lem:dM_props}, $d_{M\dual v} \geq d_M$ and so
$d_{M\dual v} = d_M$ and $Z \in \minc(M\dual v)$.

Consider now $M'' = M \dual v$. We have seen that $d_{M\dual v}
= d_M$, and so $d_{M''} = d_M < d_{M*v} = d_{M*v+v} = d_{M\dual
v * v} = d_{M''*v}$. By the first part of this proof we have
that $Z \in \minc(M'') = \minc(M\dual v)$ implies $Z \in
\minc(M'' \dual v) = \minc(M)$. Consequently, we obtain
$\minc(M\dual v) = \minc(M)$.
\end{Proof}

From this result we see that the values of $d_M$, $d_{M*v}$,
and $d_{M\dual v}$ are either (1) all equal or (2) of the form
$m$, $m$, $m+1$ (in some order). We show in
Section~\ref{sec:distance_DM} that for \dmatroids only the
latter case occurs.

\begin{Example}
Let $V = \{a,b,c\}$, and let $M$ be the set system $(V,
\{\{a\}, \{b,c\}\})$. We have $M*b = (V, \{\{a,b\}, \{c\}\})$
and $M\dual b = (V, \{\{a\}, \{c\}, \{b,c\}\})$. Hence $d_M =
d_{M*b} = d_{M\dual b} = 1$.
\end{Example}

\section{A Characterization of Delta-Matroids}
\label{sec:char_dmatroids} By Lemma~\ref{lem:dM_nullity_prop},
$d_M(X) = d_{M*X}$ is the minimal cardinality of the sets in
$M*X$. As a consequence, each set in $M*X$ of cardinality
$d_M(X)$ belongs to $\min(M*X)$, but the converse does not
necessarily hold, i.e., the inclusion $\minc(M*X) \subseteq
\min(M*X)$ may not be an equality. We consider now set systems
with the property that the converse \emph{does} hold: for each
$X \subseteq V$, the sets in $\min(M*X)$ are all of equal
cardinality (or equivalently, $\minc(M*X) = \min(M*X)$).
\begin{Definition}
A proper set system $M$ over $V$ is called \emph{isodistant} if, for
each $X \subseteq V$, $\min(M*X)$ is equicardinal.
\end{Definition}

Thus for isodistant $M$, the common cardinality of the sets in
$\min(M*X)$ is equal to $d_M(X)$. As we have noted this also holds
the other way around, and so the minimal sets in $M*X$ are
characterized by their cardinality.

Clearly, the isodistant property of set systems is invariant under
pivot: if set system $M$ is isodistant, then $M*X$ is isodistant for
each $X \subseteq V$. Due to the duality $\min(M)*V = \max(M*V)$,
one easily verifies that $M$ is isodistant iff for each $X \subseteq
V$, $\max(M*X)$ is equicardinal.
In that case the sets in $\max(M*X)$ are
all of cardinality equal to $|V| - d_M(V\setminus X)$.

A \dmatroid is a proper set system $M$ that satisfies the
\emph{symmetric exchange axiom}: For all $X,Y \in M$ and all $u
\in X \sdif Y$, either $X \sdif \{u\} \in M$ or there is a $v
\in X \sdif Y$ with $v \not= u$ such that $X \sdif \{u,v\} \in
M$ \cite{bouchet1987}. The notion of \dmatroid is equivalent to
the notion of \emph{Lagrangian matroid}
\cite[Section~6]{symplectic/matroids/borovik98}. If we assume a
matroid $M$ is described by a tuple $(V,B)$ where $B$ is the
set of bases of $M$, then it is shown in
\cite[Proposition~3]{DBLP:conf/ipco/Bouchet95} that a matroid
$M$ is precisely a equicardinal \dmatroid (the result
essentially follows from
\cite[Theorem~1]{CommentBases/Brualdi69}). It is stated in
\cite[Property~4.1]{Bouchet_1991_67} that a set system $M$ is a
\dmatroid iff $\max(M * X)$ is a matroid (described by its
bases) for every $X \subseteq V$. Consequently, every \dmatroid
is isodistant.

We now show that, surprisingly, the converse holds. Hence, the
notions of \dmatroid and isodistance are equivalent, i.e.,
without assuming the matroid structure of the maximal or minimal elements.
\begin{Theorem} \label{thm:isodistant_dmatroid_eq}
Let $M$ be a proper set system. Then $M$ is a \dmatroid iff
$M$ is isodistant.
\end{Theorem}
\begin{Proof}
Assume first that $M$ is isodistant. Let $X, Y \in M$ and $u
\in X \sdif Y$. We need to show that either $X \sdif \{u\} \in
M$ or there is a $v \in X \sdif Y$ with $v \not= u$ such that
$X \sdif \{u,v\} \in M$. Consider $M' = M*(X \sdif \{u\})$. If
$\emptyset \in M'$, then $X \sdif \{u\} \in M$ and we are done.
Assume $\emptyset \not\in M'$. We have $\{u\} \in M'$ and $Z =
Y \sdif (X \sdif \{u\}) \in M'$. As $M$ is isodistant, so is
$M'$ and thus $\{v\} \in M'$ for some $v \in Z$. As $u \not\in
Z$, $u \not= v$. Therefore, $X \sdif \{u,v\} \in M$ and we are
done.

The forward implication, i.e., the fact that the maximal
elements of a \dmatroid are of equal cardinality, follows from
\cite[Property~4.1]{Bouchet_1991_67} (stated above) or
\cite[Lemma~6]{DBLP:journals/dm/ChandrasekaranK88}.
\end{Proof}

By restricting to equicardinal set systems we obtain the
following corollary.
\begin{Corollary} \label{cor:char_matroids_bases}
Let $M$ be a proper set system. Then $M$ is a matroid described by
its bases iff both (1) $M$ is equicardinal, and (2) for each $X
\subseteq V$, $\min(M*X)$ is equicardinal.
\end{Corollary}

Although the characterization of a matroid in
Corollary~\ref{cor:char_matroids_bases} is novel, we can link
it to a well-known characterization of matroids $M$ given
below, where $M$ is described by its independent sets.
This characterization can be found, e.g., in \cite[Exercise~1.1.3]{Oxley/MatroidBook-2nd} and in \cite[Section~1.5]{Welsh/MatroidBook}.

\begin{Proposition} \label{prop:matroid_card_prop}
Let $M$ be a proper set system. Then $M$ is a matroid described by
its independent sets iff both (1) for each $Y \in M$ and $Y'
\subseteq Y$, we have $Y' \in M$, and (2) for each $X \subseteq V$,
$\max(M\sub{X})$ is equicardinal.
\end{Proposition}

The second property of Proposition~\ref{prop:matroid_card_prop} is
known as the \emph{cardinality property}.

Inspired by Proposition~\ref{prop:matroid_card_prop} and
Corollary~\ref{cor:char_matroids_bases} we obtain the following
result, which from appearance may be thought of as the ``analog'' of
Corollary~\ref{cor:char_matroids_bases} where the matroid is
described by its independent sets (it appears that there is no
obvious ``analog'' of Proposition~\ref{prop:matroid_card_prop} for
matroids described by its bases).

\begin{Corollary} \label{cor:char_matroids_indep_sets}
Let $M$ be a proper set system. Then $M$ is a matroid described by
its independent sets iff both (1) for each $Y \in M$ and $Y'
\subseteq Y$, we have $Y' \in M$, and (2) for each $X \subseteq V$,
$\min(M*X)$ is equicardinal.
\end{Corollary}
\begin{Proof}
Let $M$ be a proper set system such that condition (1) holds.
Let $X \subseteq V$, and let $Z \in \min(M*X)$. Then $Z \sdif X
\in M$. If $v \in Z \setminus X$, then $(Z \setminus \{v\})
\sdif X \subset Z \sdif X$ and $(Z \setminus \{v\}) \sdif X \in
M$ by condition (1), contradicting the minimality of $Z$.
Therefore $Z \subseteq X$. Consequently, $Z \sdif X \subseteq
X$. Hence $\min(M*X) = \min(M\sub{X}*X)$.

As $X$ is the ground set of $M\sub{X}$, we have
$\min(M\sub{X}*X) = \max(M\sub{X})*X$. Again, as $X$ is the
ground set of $M\sub{X}$, $\max(M\sub{X})*X$ is equicardinal
iff $\max(M\sub{X})$ is equicardinal. The result follows now by
Proposition~\ref{prop:matroid_card_prop}.
\end{Proof}

Note that, again, $\min(M*X)$ in
Corollary~\ref{cor:char_matroids_indep_sets} may equivalently
be replaced by $\max(M*X)$. Also note that while the second
condition of Corollary~\ref{cor:char_matroids_bases} and of
Corollary~\ref{cor:char_matroids_indep_sets} are identical,
they concern (in general) very different set systems. Indeed,
if $M$ is a matroid described by its independent sets, then
$\max(M)$ is the corresponding matroid described by its bases.

From now on, we prefer the term \dmatroid instead of the
equivalent notion of isodistant set system, as the former is
well known. However, the results in this paper do not
(directly) use the definition of \dmatroid; we use only the
property of isodistance.

\section{Distance in Delta-Matroids} \label{sec:distance_DM}
We reconsider the distance function $d_{M}$, but now restricted
to \dmatroids $M$ rather than set systems in general.

We may now characterize \dmatroids through distance and
deletion.
\begin{Theorem} \label{thm:char_dmatroid_dist_remove}
Let $M$ be a proper set system. Then $M$ is a \dmatroid iff
$d_{M*Y} = d_{M*Y\sub{X}}$ for all $X,Y
\subseteq V$ with $M*Y\sub{X}$ proper.
\end{Theorem}
\begin{Proof}
We first show the forward implication. It is easy to see that
$\min(M\sub{X}) \subseteq \min(M)$ for any set system $M$ and
$X \subseteq V$ with $M\sub{X}$ proper. Hence if $M$ is a
\dmatroid, then $d_{M} = d_{M\sub{X}}$, and similarly for $M*Y$
for all $Y \subseteq V$. To show the reverse implication,
assume $M$ is not a \dmatroid, and let $Z_1,Z_2 \in \min(M*Y)$
with $|Z_1| < |Z_2|$. We have that $M*Y\sub{Z_2}$ consists only
of $Z_2$ and therefore $d_{M*Y\sub{Z_2}} = |Z_2|$. However,
$d_{M*Y} \leq |Z_1|$ --- a contradiction.
\end{Proof}

It follows from Theorem~\ref{thm:char_dmatroid_dist_remove}
that the distance function $d_{M}$ behaves well under removal
of elements from the ground set $V$.
\begin{Theorem} \label{thm:iso_invariant_minors}
Let $M$ be a \dmatroid, and $X \subseteq V$. If $M \sub{X}$ is
proper, then $d_{M \sub{X}}(Y) = d_{M}(Y)$ for all $Y \subseteq X$.
\end{Theorem}
\begin{Proof}
We have $d_{M}(Y) = d_{M*Y}$ and $d_{M \sub{X}}(Y) =
d_{M\sub{X}*Y} = d_{M*Y\sub{X}}$ where in the last equality we
use $Y \subseteq X$ and the commutation of vertex flips. The
result holds by Theorem~\ref{thm:char_dmatroid_dist_remove}.
\end{Proof}

In particular, by Theorem~\ref{thm:iso_invariant_minors},
$d_{M}(X) = d_{M \sub{X}}(X)$ for all $X \subseteq V$ where $M
\sub{X}$ is proper, hence $d_M(X) = \min(\{|X \setminus Y| \mid
Y \in M, Y \subseteq X \})$.

The property of \dmatroids shown in
Theorem~\ref{thm:iso_invariant_minors} is important in a
subsequent study of interlace polynomials on \dmatroids
\cite{BH/InterlacePolyDM/10}. Of course,
Theorem~\ref{thm:iso_invariant_minors} does not hold for set
systems in general. Indeed, it is easy to verify that set
system $M = (V,\{\emptyset,V\})$ is not a \dmatroid for $|V|
\geq 3$. Take $|V|=3$. We have, for $u \in V$, $M \setminus u =
(V\setminus \{u\},\{\emptyset\})$ and therefore $2 = d_{M
\setminus u}(V\setminus \{u\}) \neq d_{M}(V\setminus \{u\}) =
1$. It is also easy to verify that the property of
Theorem~\ref{thm:iso_invariant_minors} does not characterize
\dmatroids like in Theorem~\ref{thm:char_dmatroid_dist_remove}
(take, e.g., $M = (\{a,b,c\},\{\emptyset,\{a\},\{b,c\}\})$).

By Theorem~\ref{thm:eq_m_dualp_setsystem}, for arbitrary set
systems the value of $d_{M\dual v}$ is the minimum of $d_M$ and
$d_{M*v}$ when these two values differ. However, the value of
$d_{M\dual v}$ could not be fixed when $d_M$ equals $d_{M*v}$.
This changes when $M$ is a \dmatroid.
\begin{Lemma} \label{lem:min_vxor_equal_nullity}
Let $M$ be a \dmatroid, and $v \in V$ such that $d_M =
d_{M*v}$. Then {\em (1)} no set in $\min(M)$ contains $v$, {\em
(2)} $\min(M) = \min(M*v)$, and {\em (3)} $d_{M\dual v} =
d_M+1$.
\end{Lemma}
\begin{Proof}
Let $m = d_M = d_{M*v}$.
(1) Let $Y \in \min(M)$. Then $Y \sdif\{v\} \in M*v$. As $|Y|
\leq |Y \sdif\{v\}|$, we have $v \notin Y$.

(2) Let $Y \in \min(M)$. As $Y \cup \{v\} \in M*v$, there must
be a $Y' \in \min(M*v)$ of cardinality $m$ with $Y'\subseteq
Y$. If $Y'= Y\cup\{v\} \setminus \{w\}$ with $v\neq w$, then
$Y'\setminus \{w\} \in \min(M)$ while this set is smaller than
$m$ --- a contradiction. Hence $v=w$, and $Y = Y' \in
\min(M*v)$. Therefore, $\min(M) \subseteq \min(M*v)$. By
symmetry we obtain the other inclusion.

(3) By Lemma~\ref{lem:dM_props}, $d_{M\bar*v} \in \{m,m+1\}$,
so it suffices to prove that there are no sets in $M\dual v$
that have cardinality $m$. Thus assume $Z \in M\dual v$ and
$|Z| = m$. By the definition of dual pivot, either $Z\in M$
(case $v\in Z$) or $Z$ in exactly one of $M$ and $M*v$ (case
$v\notin Z$). In the former case we have $Z \in \min(M)$ while
$v \in Z$, contradicting (1). The latter case contradicts with
(2).
\end{Proof}

It is observed in \cite{BH/PivotLoopCompl/09} that both loop
complementation and dual pivot do not (in general) retain the
property of being a \dmatroid. For example, for \dmatroid $M =
(V,2^V \setminus \{V\})$ with $V = \{1,2,3\}$, $M\dual 1 $ is
not a \dmatroid. In fact, $\min(M \dual 1) =
(V,\{\{1\},\{2,3\}\})$ is not even equicardinal. The next
result shows that $\min(M \dual v)$ is equicardinal for a
\dmatroid $M$ when $d_M \not= d_{M*v}$.

\begin{Lemma} \label{lem:dmat_dual_min_equi}
Let $M$ be a \dmatroid, and $v \in V$. If $d_M \not= d_{M*v}$,
then $\min(M\dual v)$ is equicardinal.
\end{Lemma}
\begin{Proof}
Assume first that $d_M < d_{M*v}$. By
Theorem~\ref{thm:eq_m_dualp_setsystem}, $d_M = d_{M\dual v}$,
and $\min(M) = \minc(M\dual v) \subseteq \min( M\dual v )$. Let
$Z \in \min(M \dual v)$. By definition of dual pivot, $Z \in M$
or $Z\sdif \{v\} \in M$ (or both). Hence, as no set in
$\min(M)$ contains $v$, there is a $Y \in \min(M)$ with $Y
\subseteq Z \setminus \{v\}$. Again, $Y \in \min(M)$ implies $Y
\in \min(M \dual v)$ and so $Y = Z$.

The other case, $d_M > d_{M*v}$, follows by symmetry (consider
the \dmatroid $M' = M*v$ similar as in the proof of
Theorem~\ref{thm:eq_m_dualp_setsystem}).
\end{Proof}

Let $M = (V,D)$ be a set system and $v \in V$. We define
\emph{pseudo-deletion} of $M$ on $v$, denoted by $M
\widehat\vertexrem v$, as $M \widehat\vertexrem v
=(V,D\vertexrem v*v)$. Similarly, we define
\emph{pseudo-contraction} of $M$ on $v$, denoted by $M
\widehat\slash v$, as $(V,D*v\vertexrem v)$. Note that the
ground sets for both pseudo-deletion and pseudo-contraction
remain unchanged. Also note that $M*v\widehat\vertexrem v =
M\widehat\slash v *v$.

The definitions of pseudo-deletion and pseudo-contraction are
motivated by matroids as follows. Recall that for a matroid $M$
described by its bases and $v \in V$, $M\vertexrem v$ and
$M*v\vertexrem v$ are the matroid operations of \emph{deletion}
(if $v$ is not a coloop) and \emph{contraction} (if $v$ is not
a loop), denoted by $M\vertexrem v$ and $M \slash v$,
respectively. It is easy to see that then pseudo-deletion is
adding $v$ as a coloop to $M\vertexrem v$ and
pseudo-contraction is adding $v$ as a loop to $M\slash v$
\cite{Bouchet/87/ejc/isotropicsys}. In this way, we regard
pseudo-deletion and pseudo-contraction as matroid operations as
well. Pseudo-deletion and pseudo-contraction take the following
form
if a matroid is described by its circuits. If $M' =
(V,\mathcal{C})$ is the circuit description of $M$, then
$M'\widehat\vertexrem v = (V,\mathcal{C}\vertexrem v)$ and
$M'\widehat\slash v = (V,(\mathcal{C}\slash v) \cup
\{\{v\}\})$.

We are now ready to formulate the announced $m,m,m+1$ result for \dmatroids.
\begin{Theorem} \label{thm:null_plus_op_setsystem}
Let $M$ be a \dmatroid, and $v \in V$. Then the equicardinal
set systems $\min(M)$, $\min(M*v)$, and $\minc(M\dual v)$ are
such that precisely two of the three are equal, to say $M_1$.
Moreover, the third $M_2$ is such that $d_{M_2} = d_{M_1}+1$
and $M_2 \widehat\slash v = M_1$.

In particular, the values of $d_{M}$, $d_{M*v}$, and $d_{M\dual v}$
are such that precisely two of the three are equal, to say $m$, and
the third is equal to $m+1$.
\end{Theorem}
\begin{Proof}
(i) The case $d_{M} = d_{M*v}$ follows from
Lemma~\ref{lem:min_vxor_equal_nullity} except for the equality
$M_2 \widehat\slash v = M_1$. Let $Z \in M_1 = \min(M*v)$. By
Lemma~\ref{lem:min_vxor_equal_nullity}, $v \not\in Z$. We have
$Z\sdif\{v\} \in M$ and $v \in Z\sdif\{v\}$, and therefore $Z
\sdif \{v\} \in M \dual v$ and $|Z \sdif \{v\}| = d_{M \dual
v}$. Hence, $Z \sdif \{v\} \in \minc(M \dual v) = M_2$.
Conversely, let $Z \in \minc(M \dual v) = M_2$ with $v \in Z$.
Then $Z \in M$. Hence $Z \setminus \{v\} \in M*v$. Since
$d_{M_2} = d_{M_1}+1$, $Z \setminus \{v\} \in \min(M*v) = M_1$.
(ii) Consider $d_{M} < d_{M*v}$, hence $d_{M}+1 = d_{M*v}$. By
Theorem~\ref{thm:eq_m_dualp_setsystem} we know that
$\minc(M\dual v) = \min(M)$. If $Z \in \min(M)$, then $Z\sdif v
\in M*v$ is minimal by cardinality and so $\min(M)*v \subseteq
\min(M*v)$. Conversely, if $Z \in \min(M*v)$ and $v \in Z$,
then $Z \setminus \{v\} \in \min(M)$ as $d_{M_2} = d_{M_1}+1$.
(iii) The case $d_{M*v} < d_{M}$ holds by symmetry.
\end{Proof}

Theorem~\ref{thm:null_plus_op_setsystem} is related to
Theorem~(9.4) of \cite{Bouchet/87/ejc/isotropicsys}, which
deals fundamentally with binary matroids. In fact, Theorem~\ref{thm:null_plus_op_setsystem}
may be seen as a generalisation of Theorem~(9.4) of \cite{Bouchet/87/ejc/isotropicsys},
cf. Theorem~\ref{thm:nullspace_plus_op_graph}.

Note that the cardinality of the sets in $M_2$ is exactly one
larger than the cardinality of the sets in $M_1$. Also note
that $\min(M)$ and $\min(M*v)$ are matroids and with ranks
$d_{M}$ and $d_{M*v}$, respectively.
Lemma~\ref{lem:dmat_dual_min_equi} shows that if $d_M \not=
d_{M*v}$, then $\minc(M\dual v) = \min(M\dual v)$ is also a
matroid with rank $d_{M\dual v}$. Of course, if $M\dual v$ is a
\dmatroid, then $\minc(M\dual v) = \min(M\dual v)$ also holds.

Let $X\subseteq V$. By
Theorem~\ref{thm:null_plus_op_setsystem}, if $v \not\in X$,
then the values of $d_{M}(X)$, $d_{M*v}(X)$, and
$d_{M\dual v}(X)$ are such that precisely two of the three are
equal, to say $m$, and the third is equal to $m+1$. Also, if $v
\in X$, then the same statement holds for $d_{M}(X)$, $d_{M*v}(X)$, and
$d_{M+v}(X)$.

We may state the $\max$ analog of
Theorem~\ref{thm:null_plus_op_setsystem},
using the duality $\max(M)*V = \min(M*V)$ --- note again the
change from $M\dual v$ to $M+v$.
\begin{Theorem} \label{thm:nullspace_plus_op_setsystem_max}
Let $M$ be a \dmatroid, and $v \in V$. Then the equicardinal
set systems $\max(M)$, $\max(M*v)$, and $\maxc(M+v)$ are such
that precisely two of the three are equal, to say $M_1$.
Moreover, the third $M_2$ is such that $d_{M_1} = d_{M_2}+1$
and $M_2 \widehat\vertexrem v = M_1$.
\end{Theorem}
\begin{Proof}
We have $\max(M) = \min(M*V)*V$, $\max(M*v) = \min((M*V)*v)*V$,
and $\max(M+v) = \min(M+v*V)*V = \min((M*V)\dual v)*V$. Now
apply Theorem~\ref{thm:null_plus_op_setsystem} to \dmatroid
$M*V$. Finally, let $M_1 = M'_1*V$ and $M_2 = M'_2*V$, where
$M'_1$ and $M'_2$ are the two set systems of
Theorem~\ref{thm:null_plus_op_setsystem} belonging to \dmatroid
$M*V$. We have $M_2 \widehat\vertexrem v = M'_2 * V
\widehat\vertexrem v = M'_2 \widehat\slash v *V = M'_1 * V =
M_1$.
\end{Proof}

Of course, $\max(M)$ and $\max(M*v)$ are matroids. If $M+v$ is
a \dmatroid, then $\max(M+v) = \maxc(M+v)$ is also a matroid.
The matroid formulation of
Theorem~\ref{thm:nullspace_plus_op_setsystem_max} for the case
where $M+v$ is a \dmatroid is given in the Introduction.

Note that $d_{M_1} = d_{M_2}+1$ means that the cardinality of
the sets in $M_2$ is exactly one smaller than the cardinality
of the sets in $M_1$.

Note also that the definition of loop complementation may be
formulated through the distance measure since $d_{M+v}(X)=0$
iff $X \in M+v$. We have therefore, for $X \subseteq V$, by
definition of loop complementation, $d_{M+v}(X)=0$ iff
$d_{M}(X)=0$ when $v \not\in X$, and $d_{M+v}(X)=0$ iff
$(d_{M}(X)=0) \xor (d_{M*v}(X)=0)$ when $v \in X$. Recall that
$d_{M+v}(X) = d_{M}(X)$, hence the case $v \not\in X$ is
extended to arbitrary values of $d_{M+v}(X)$. Moreover, by
Theorem~\ref{thm:nullspace_plus_op_setsystem_max}, the case $v
\in X$ is extended for \dmatroids to arbitrary values of
$d_{M+v}(X)$ through the $m$, $m$, $m+1$ property (by extending
$\xor$ in a suitable way from Booleans to integers).

One may wonder whether or not the property of
Theorem~\ref{thm:null_plus_op_setsystem} characterizes \dmatroids.
The next example illustrates that this is not the case.
\begin{Example} \label{ex:no_char_dmatroids_kkko}
Let $M = (V,\{\emptyset,V\})$. Recall that $M$ is not a
\dmatroid for $|V| \geq 3$. Assume now that $|V| \geq 3$ is
even. Let $X \subseteq V$ and $v \in X$. As the (two) sets in
$M$ are of equal parity, the distances $d_{M}(X)$ and
$d_{M*v}(X) = d_{M}(X \sdif \{v\})$ are of different parity (as
$|X|$ and $|X \sdif \{v\}|$ are of different parity). Now, by
Lemma~\ref{lem:dM_props}, $d_{M+v}(X) =
\min\{d_{M}(X),d_{M*v}(X)\}$, and we have that $d_{M}(X)$,
$d_{M*v}(X)$, and $d_{M+v}(X)$ are, in this order, either of
the form $m$, $m+1$, and $m$, or of the form $m+1$, $m$, and
$m$.
\end{Example}

\section{Representable Delta-Matroids} \label{sec:repres_sets}
In this section we consider the case where a \dmatroid $M$ is
represented by a matrix $A$. We show that in that case the
notion of distance to $X$ in the represented \dmatroid closely
matches that of nullity of the principal submatrix $A[X]$.

For a $V\times V$ matrix $A$ (the columns and rows of $A$ are
indexed by finite set $V$) and $X \subseteq V$, $A[X]$ denotes
the principal submatrix of $A$ with respect to $X$, i.e., the $X \times
X$ matrix obtained from $A$ by restricting to rows and columns
in $X$. We also define $A \setminor X = A \sub{V \setminus X}$.
For $V \times V$ matrix $A$ we consider the associated set
system $\mathcal{M}_A = (V,D_A)$ with $D_A = \{ X \subseteq V
\mid A\sub{X} \mbox{ is nonsingular}\}$. Observe that
$\mathcal{M}_{A\sub{X}} = \mathcal{M}_A\sub{X}$, and
$\mathcal{M}_{A\setminus X} = \mathcal{M}_{A} \setminus X$.

It is shown in \cite{bouchet1987} that $\mathcal{M}_A$ is a
\dmatroid when $A$ is symmetric or skew-symmetric (over some
field $\mathbb{F}$). Note that $\varnothing\in
\mathcal{M}_{A}$. We say that \dmatroid $M$ is
\emph{representable over $\mathbb{F}$} if $M = \mathcal{M}_A*
X$ for some skew-symmetric matrix $A$ and $X \subseteq V$; $A$
is called a \emph{representation} of $M$. A \dmatroid is called
\emph{binary} if it is representable over $\two$.

Recall that for a $W \times V$ matrix $A$, the \emph{column
matroid} $N = (V,B)$ of $A$ described by its bases is such
that, for $X \subseteq V$, $X \in B$ iff the columns of $A$
belonging to $X$ form a basis of the column space of $A$.
Matrix $A$ is said to \emph{represent} $N$. For matroids, this
usual sense of representability coincides with representability
in the \dmatroid sense, see 4.4 of \cite{bouchet1987}. Hence,
every binary matroid is a binary \dmatroid.

We now formulate the matroid version of the strong principal
minor theorem \cite{Kodiyalam_Lam_Swan_2008} (the original
result is more general, as it considers quasi-symmetric
matrices over a division ring), see also
\cite[Lemma~10]{MaxPivotsGraphs/Brijder09}.
\begin{Proposition}[Strong Principal Minor Theorem \cite{Kodiyalam_Lam_Swan_2008}] \label{prop:spmt}
Let $A$ be a $V \times V$ symmetric or skew-symmetric matrix
(over some field). Then $\max(\mathcal{M}_A)$ equals the
column matroid of $A$ (described by its bases).
\end{Proposition}

As a consequence of the strong principal minor theorem, the
sets in $\max(\mathcal{M}_A)$ are all of cardinality equal to
the rank $r(A)$ of $A$ --- this fact is known as the principal
minor theorem.
We now use the principal minor theorem to obtain that the
distance $X \subseteq V$ to $\mathcal{M}_A$ corresponds to the
nullity of $A\sub{X}$. We will also use the strong principal minor
theorem to prove Theorem~\ref{thm:nullspace_plus_op_graph}.
\begin{Theorem} \label{thm:respr_setsystem_isodistant_skew_symm}
Let $A$ be a $V \times V$ symmetric or skew-symmetric matrix
(over some field). Then
$d_{\mathcal{M}_A}(X) = n(A\sub{X})$ for every $X \subseteq V$.
\end{Theorem}
\begin{Proof}
By Theorem~\ref{thm:iso_invariant_minors},
$d_{\mathcal{M}_A}(X) = d_{\mathcal{M}_A\sub{X}}(X)$ (note that
$\mathcal{M}_A\sub{X}$ is proper as $\emptyset \in
\mathcal{M}_A$). Now, $d_{\mathcal{M}_A\sub{X}}(X) = |X| - |Z|$
with $Z \in \max(\mathcal{M}_A\sub{X})$. Moreover,
$\mathcal{M}_A\sub{X} = \mathcal{M}_{A\sub{X}}$. By the
principal minor theorem, $|Z| = r(\mathcal{M}_{A\sub{X}}) =
r(A\sub{X})$, and so $d_{\mathcal{M}_{A\sub{X}}}(X) =
n(A\sub{X})$.
\end{Proof}

To extend the notion of nullity to \dmatroids (or proper set
systems in general), we regard $d_M(X)$ as the \emph{nullity}
of $X$ in $M$. We may now also define the \emph{rank} $r_M(X)$
of $X$ in $M$ by $r_M(X) = |X| - d_M(X)$. In this way we have
$r_{\mathcal{M}_A}(X) = r(A\sub{X})$, where $r(A\sub{X})$
denotes the rank of matrix $A\sub{X}$.

\begin{Remark}
There have been a number of rank functions introduced for
\dmatroids. In \cite{DBLP:journals/combinatorics/BouchetJ00}
the rank of $X \subseteq V$ is defined as $r'_M(X) = \max\{ |X
\cap Y| + |(V \setminus X) \cap (V \setminus Y)| \mid Y \in
M\}$. It is easy to verify that $r'_M(X) = \max \{ |V \setminus
(X \sdif Y)| \mid Y \in M \} = |V| - \min\{ |X \sdif Y| \mid Y
\in M \}$. Therefore, $r'_M(X) = |V| - d_M(X)$ and hence the
notion is slightly different from the rank function $r_M(X)$
defined in this paper. Also, in \cite{Bouchet1989/maps_deltam}
the birank of $X,Y \subseteq V$ with $X \cap Y = \emptyset$ is
defined as $r''_M(X,Y) = \max \{ |Z \cap X| + |(V \setminus Z)
\cap Y| \mid Z \in M \}$. We have $r''_M(X,V \setminus X) =
r'_M(X)$. Finally, in \cite{DBLP:journals/jamds/KabadiS05} the
birank of $X,Y \subseteq V$ with $X \cap Y = \emptyset$ is
defined as $r'''_M(X,Y) = \max \{ |Z \cap X| - |Z \cap Y| \mid
Z \in M \}$. Function $r'''_M(X,Y)$ does correspond to the rank
$r_M(X)$ as defined in this paper, as $r'''_M(X,V \setminus X)
= r_M(X)$. However none of these papers on \dmatroids
(explicitly) considers nullity as a distance measure.
\end{Remark}

A \dmatroid $M$ is called \emph{even} if the cardinality
of the sets in $M$ all have equal parity. Let $A$ be a
skew-symmetric matrix over some field $\mathbb{F}$. It easily follows from \cite[Thm~4.3.3]{bouchet1987} that \dmatroid $\mathcal{M}_A$ is even when $A$ is zero-diagonal
(note that this condition is only relevant when $\mathbb{F}$ is of characteristic $2$).
We now obtain the following corollary to Theorem~\ref{thm:respr_setsystem_isodistant_skew_symm}.

\begin{Corollary} \label{cor:zerodiag_skew_nullity}
Let $A$ be a $V \times V$ zero-diagonal skew-symmetric matrix, and $v \in V$. Then $n(A)$ and $n(A \vertexrem v)$ differ by precisely $1$.
\end{Corollary}
\begin{Proof}
Since the cardinality of the sets in $\mathcal{M}_A$ have a
common parity, for all $X \subseteq V$, $d_{\mathcal{M}_A}(X)$
is odd iff $d_{\mathcal{M}_A}(X\sdif\{v\})$ is even. As
$|d_{\mathcal{M}_A}(X) - d_{\mathcal{M}_A}(X\sdif\{v\})| \leq
1$, we have $|d_{\mathcal{M}_A}(X) -
d_{\mathcal{M}_A}(X\sdif\{v\})| = 1$. Let now $X = V$ and we
obtain by Theorem~\ref{thm:respr_setsystem_isodistant_skew_symm}, $|n(A) - n(A \vertexrem v)| = 1$.
\end{Proof}

\section{Application: Graphs} \label{sec:null_graphs}
In this section we translate the above results to the realm of
graphs, where, e.g., the operations ${}+v$ and ${}*v$ have
their own specific interpretation. We consider undirected
graphs without parallel edges, but we do allow loops. For a
graph $G = (V,E)$ and $x \in V$, we have $\{x\} \in E$ iff $x$
has a loop. With a graph $G$ one associates its adjacency
matrix $A(G)$, which is a $V \times V$ matrix
$\left(a_{u,v}\right)$ over $\two$ with $a_{u,v} = 1$ iff
$\{u,v\} \in E$ (with possibly $u=v$). In this way, the family
of graphs with vertex set $V$ corresponds precisely to the
family of symmetric $V \times V$ matrices over $\two$.
Therefore we often make no distinction between a graph and its
matrix, so, e.g., by the null space (or kernel) and nullity
(i.e., dimension of the null space) of graph $G$, denoted by
$\ker(G)$ and $n(G)$ respectively, we mean the null space and
nullity of its adjacency matrix $A(G)$ (computed over $\two$).
Also, for $X \subseteq V$, $G\sub{X} = A(G)\sub{X}$ is the
subgraph of $G$ induced by $X$. By convention, the empty
graph/matrix is nonsingular. Similar as for set systems, we
often write $V$ to denote the vertex set of the graph under
consideration.

For a graph $G$ and a set $X \subseteq V$, the graph obtained
after \emph{loop complementation} for $X$ on $G$, denoted by
$G+X$, is obtained from $G$ by adding loops to vertices $v \in
X$ when $v$ does not have a loop in $G$, and by removing loops
from vertices $v \in X$ when $v$ has a loop in $G$. Hence, if
one considers a graph as a matrix, then $G+X$ is obtained from
$G$ by adding the $V \times V$ matrix with elements $x_{i,j}$
such that $x_{i,j} = 1$ if $i=j\in X$ and $0$ otherwise. Note
that $(G+X)+Y = G+(X \sdif Y)$.

Given the set system $\mathcal{M}_G = \mathcal{M}_{A(G)} =
(V,D_G)$ for some graph $G = (V,E)$, one can (re)construct the
graph $G$, see \cite[Property~3.1]{Bouchet_1991_67}. Hence the
function $\mathcal{M}_{(\cdot)}$ which assigns to each graph
$G$ its set system $\mathcal{M}_G$ is injective. In this way,
the family of graphs (with set $V$ of vertices) can be
considered as a subset of the family of set systems (over set
$V$). Note that $\mathcal{M}_{(\cdot)}$ is not injective for
matrices over $\two$ in general: e.g., for fixed $V$ with $|V|
= 2$, the $2 \times 2$ zero matrix and the matrix
$\left(\begin{array}{cc} 0 & 1 \\
0 & 0 \end{array}\right)$ correspond to the same set system.

It is shown in \cite{BH/PivotLoopCompl/09} that
$\mathcal{M}_{G+X} = \mathcal{M}_{G}+X$ for any graph $G$ and
$X \subseteq V$. Therefore the operation ${}+ X$ on set systems
$M$ is a generalization of loop complementation on graphs $G$
---
which explains its name.

If $G$ is a graph and $u$ a vertex\footnote{Local
complementation is often defined on simple graphs; here we
consider the obvious extension to graphs where loops are
allowed. Note that local complementation may be applied here to
a non-looped vertex $u$, which is different from, e.g.,
\cite{BH/PivotLoopCompl/09}.} of $G$, then the result of
\emph{local complementation} of $u$ on $G$, denoted by
$\mathrm{loc}_u(G)$, is the graph obtained from $G$ by
``toggling'' the edges in the neighbourhood $N_G(u) = \{ v \in
V \mid \{u,v\} \in E(G), u \not= v \}$ of $u$ in $G$: for each
$v,w \in N_G(u)$, $\{v,w\}\in E(G)$ iff $\{v,w\} \not\in
E(\mathrm{loc}_u(G))$ (again, $v=w$ is possible). The other
edges are left unchanged.

If $u$ is a looped vertex of $G$, then it is shown in
\cite{Geelen97} that $\mathcal{M}_{G} * u =
\mathcal{M}_{\mathrm{loc}_u(G)}$. Moreover, if $u$ is a
unlooped vertex of $G$, then $\mathcal{M}_{G} \dual u =
\mathcal{M}_{\mathrm{loc}_u(G)}$ (see
\cite{BH/PivotLoopCompl/09}). In this way, local
complementation is defined for \dmatroids. For convenience, we
define the \emph{pivot} of a looped vertex $u$ on $G$, denoted
as $G*u$, by $\mathrm{loc}_u(G)$ (it is not defined on unlooped
vertices). Similarly, we define the \emph{dual pivot} of an
unlooped vertex $u$ on $G$, denoted as $G \dual  u$, by
$\mathrm{loc}_u(G)$ (it is not defined on looped vertices).
Thus, if $u$ is looped, then $\mathcal{M}_{G}
* u = \mathcal{M}_{G*u}$, and if $u$ is
unlooped, then $\mathcal{M}_{G} \dual u = \mathcal{M}_{G \dual
u}$. In general, for a set $X \in \mathcal{M}_{G}$,
$\mathcal{M}_{G} * X = \mathcal{M}_{G*X}$ where $G*X$ is a
graph called the \emph{pivot} (or \emph{principal pivot
transform}) \cite{Tsatsomeros2000151} of $G$ on $X$
\cite{bouchet1987} --- clearly, $\mathcal{M}_{G} * X$ does not
correspond to a graph if $X \not\in \mathcal{M}_{G}$.

As usual, a vector $v$ indexed by $V$ over $\two$ may be
identified by a subset $X \subseteq V$, where $v$ at position
$x \in V$ is $1$ iff $x \in X$. In this way, we regard
$\ker(G)$ for a graph $G$ as a subset of $2^V$. By the strong
principal minor theorem we have that the cycle space of the
matroid $\max(\mathcal{M}_G)$ is precisely $\ker(G)$. Of
course, the nullity of the matroid $\max(\mathcal{M}_G)$ is the
dimension $\dim(\ker(G))$ of $\ker(G)$.

If we restrict now
Theorem~\ref{thm:nullspace_plus_op_setsystem_max} for the case
where the \dmatroid $M$ is equal to $\mathcal{M}_{G}$, then we
obtain the following result.

\begin{Theorem} \label{thm:nullspace_plus_op_graph}
Let $G$ be a graph having a looped vertex $v$. Then $\ker(G)$,
$\ker(G*v)$, and $\ker(G+v)$ are such that precisely two of the
three are equal, to say $K_1$, and the third, $K_2$, is such
that $\dim(K_2) = \dim(K_1) + 1$ and $K_1 = \{X \in K_2 \mid v
\not\in X\}$.
\end{Theorem}
\begin{Proof}
From Theorem~\ref{thm:respr_setsystem_isodistant_skew_symm} we
know that for any graph $G'$ the dimension of the kernel equals
$n(G') = d_{\mathcal{M}_{G'}}(V)$, i.e., the cardinality of
sets in $\min(\mathcal{M}_{G'}*V)$ which are complements of
sets in $\max(\mathcal{M}_{G'})$.

By Theorem~\ref{thm:nullspace_plus_op_setsystem_max},
$\max(\mathcal{M}_G)$, $\max(\mathcal{M}_{G*v})$, and
$\max(\mathcal{M}_{G+v})$ are such that precisely two of the
three are equal, to say $M_1$, and the nullity of the third,
$M_2$, is one larger than the nullity of $M_1$. Moreover, the
family of circuits of $M_1$ is obtained from the family of
circuits of $M_2$ by removing the sets containing $v$.

Hence (by discussion above), $\ker(G)$, $\ker(G*v)$, and
$\ker(G+v)$ are such that precisely two of the three are equal,
to say $K_1$. The third, $K_2$, is such that $\dim(K_2) =
\dim(K_1) + 1$ and $K_1 = \{X \in K_2 \mid v \not\in X\}$.
\end{Proof}

Theorem~\ref{thm:nullspace_plus_op_graph} is similar to a
result of Traldi~\cite[Lemma~23]{LT/LAA/2011}, where graph
$G'$, obtained from $G$ by removing all edges incident to $v$
except for the loop on $v$, is considered instead of $G*v$.
Moreover, Theorem~\ref{thm:nullspace_plus_op_graph} is
essentially \cite[Theorem~(9.4)]{Bouchet/87/ejc/isotropicsys}
for the case where $G$ is a fundamental graph of an isotropic system.

Let $G$ be a graph with looped vertex $v$. By
Theorem~\ref{thm:nullspace_plus_op_graph}, the values of
$n(G)$, $n(G*v)$, and $n(G+v)$ are such that precisely two of
the three are equal, to say $m$, and the third is equal to
$m+1$. It is shown in \cite{MaxPivotsGraphs/Brijder09} that the
adjacency matrix of $G*v\vertexrem v$ is the Schur complement
of $v$ on the adjacency matrix of $G$, and moreover it is well
known, see e.g. \cite{SchurBook2005}, that the Schur complement
retains the nullity, i.e., $n(G*v\vertexrem v) = n(G)$. Hence,
we have $n(G*v) = n(G\vertexrem v)$ and we obtain as a
consequence of Theorem~\ref{thm:nullspace_plus_op_graph} the
following result of \cite{DBLP:journals/ejc/BalisterBCP02}.

\begin{Proposition}[Lemma~2 of \cite{DBLP:journals/ejc/BalisterBCP02}] \label{prop:balister_null_plus_op_graph}
Let $G$ be a graph and $v \in V$. Then the values of $n(G)$, $n(G
\vertexrem v)$, and $n(G+v)$ are such that precisely two of the
three are equal, to say $m$, and the third is equal to $m+1$.
\end{Proposition}


Corollary~\ref{cor:zerodiag_skew_nullity} for the case where the field is $\two$ may be stated as follows.
\begin{Corollary}
Let $G$ be a graph, and $v \in V$ a vertex of $G$. If $G$ has no loops,
then $n(G)$ and $n(G \vertexrem v)$ differ by precisely $1$.
\end{Corollary}

\section{Vertex-Flip-Safe Delta-Matroids} \label{sec:strong_DM}
Recall from Section~\ref{sec:distance_DM} that the result of
applying dual pivot or loop complementation on a \dmatroid is
not necessarily a \dmatroid. In this section we consider
families of \dmatroids that are closed under invertible vertex
flips in general. In particular, we show that the binary \dmatroids
form one such family. As a consequence, for binary \dmatroid
$M$ the set systems $M*v$ and $M\dual v$ in
Theorem~\ref{thm:null_plus_op_setsystem} are binary \dmatroids.

\begin{Definition} \label{def:vfclosed}
Let $M$ be a \dmatroid. We say that $M$ is a
\emph{vertex-flip-safe} (or \emph{\vfsafe } for short) if for
any sequence $\varphi$ of invertible vertex flips
(equivalently, pivots and loop complementations) over $V$ we
have that $M\varphi$ is a \dmatroid.
\end{Definition}

Hence, $M$ is a \vfsafe \dmatroid iff each set system in the
orbit of $M$ under pivot and loop complementation is a
\dmatroid.

We say that a family of \dmatroids is \emph{vf-closed} if the
family is closed under invertible vertex flips. We now
show that the family of binary \dmatroids is vf-closed. First
we remark that this is not trivial. While we know that for a
graph $G$ and $X \subseteq V$, (1) $\mathcal{M}_G*X$ is a
binary \dmatroid by definition, and (2) $\mathcal{M}_G + X =
\mathcal{M}_{G+X}$ corresponds to a graph, it is not
immediately clear that, e.g., $\mathcal{M}_G*X + Y$ is a
\dmatroid for all $Y \subseteq V$ (recall that
$\mathcal{M}_G*X$ does not correspond to a graph when $X
\not\in \mathcal{M}_G$).

\begin{Theorem} \label{thm:graph_ss_strongly_isodistant}
The family of binary \dmatroids is vf-closed. In particular,
every binary \dmatroid is \vfsafe.
\end{Theorem}
\begin{Proof}
Let $M$ be a binary \dmatroid. Hence $M$ is of the form
$\mathcal{M}_G*X$ for some graph $G$ and $X \subseteq V$. Let
$\varphi$ be a sequence of invertible vertex flips over $V$.
Let $W \in \mathcal{M}_G*X \varphi$, and consider now $\varphi'
= *X\varphi*W$. By \cite[Corollary 15]{BH/PivotLoopCompl/09},
$\varphi'$ can be put in the following normal form:
$\mathcal{M}_G \varphi' = \mathcal{M}_G+Z_1*Z_2+Z_3$ for some
$Z_1,Z_2,Z_3 \subseteq V$ with $Z_1 \subseteq Z_2$. We have
$\mathcal{M}_G + Z_1 = \mathcal{M}_{G+Z_1}$. Thus
$\mathcal{M}_G+Z_1*Z_2+Z_3 = \mathcal{M}_{G+Z_1}*Z_2+Z_3$. By
construction $\emptyset \in \mathcal{M}_G \varphi'$. Hence we
have $\emptyset \in \mathcal{M}_{G+Z_1}*Z_2$. Therefore $Z_2
\in \mathcal{M}_{G+Z_1}$ and so $G+Z_1*Z_2$ is defined.
Consequently, $G' = G+Z_1*Z_2+Z_3$ is defined and
$\mathcal{M}_G \varphi' = \mathcal{M}_{G'}$. Hence $M\varphi =
\mathcal{M}_G*X\varphi = \mathcal{M}_{G'}*W$ and thus $G'$
represents $M\varphi$. Consequently, $M\varphi$ is a binary
\dmatroid.
\end{Proof}

We consider some specific matroids to illustrate the scope of this notion.
By Theorem~\ref{thm:graph_ss_strongly_isodistant}, every binary
matroid is \vfsafe. Not every matroid is a \vfsafe matroid. The
$6$-point line, i.e., $U_{2,6} = (V,\{\{u,v\} \mid u,v \in V, u
\not= v\})$ with $|V| = 6$, is not \vfsafe. Recall that $X \in
M+V$ iff the number of sets in $M[X]$ is odd. We have $V \in
U_{2,6}+V$ as the number of sets in $U_{2,6}$ is ${6 \choose 2}
= 15$ (odd), while the sets of cardinality $4$ and $5$ are not
in $U_{2,6}+V$ as ${5 \choose 2} = 10$ and ${4 \choose 2} = 6$
are even. Consequently, the symmetric exchange axiom does not
hold for $V \in U_{2,6}+V$ (as neither $V \setminus \{u\} \in
U_{2,6}+V$ nor $V \setminus \{u,v\} \in U_{2,6}+V$ for any $u,v
\in V$).

Based on Theorem~\ref{thm:isodistant_dmatroid_eq} one easily
verifies (by computer) that several small (non-binary) matroids
are \vfsafe. Such examples include the matroids $U_{2,4}$, $U_{2,5}$,
$U_{3,6}$, $Q_6$, ${\cal W}^3$, $P_8$, $P^=_8$, and Pappus. For
information on these matroids, see the Appendix on interesting
matroids in \cite{Oxley/MatroidBook-2nd}.

We turn to minors. A \emph{minor} of a \dmatroid $M$ is a proper set
system obtained from $M$ by any sequence of ${}\vertexrem v$ (deletion)
and ${}* v \vertexrem v$ (contraction) operations. A minor of $M$ is thus a
proper set system of the form $M * X \vertexrem Y$ with $X\subseteq Y
\subseteq V$. Consequently, a minor of a \dmatroid is again a \dmatroid.
Also note that this notion of minor restricted to matroids coincides with the usual
notion of minor for matroids.

\begin{Theorem} \label{thm:vfclosed_del}
The family of \vfsafe \dmatroids is minor-closed. In
particular, the family of \vfsafe matroids is minor-closed.
\end{Theorem}
\begin{Proof}
It suffices to consider $M \vertexrem u$ for a \vfsafe
\dmatroid and some $u\in V$. Let $\varphi$ be a sequence of
invertible vertex flips on $V\setminus\{u\}$. Then $M' = (M
\vertexrem u) \varphi = (M \varphi) \vertexrem u$. Moreover, $M
\varphi$ is a \dmatroid as $M$ is a \dmatroid. Also, $M'$ is
proper, as $M \vertexrem u$ is proper. Consequently, $M' = (M
\varphi) \vertexrem u$ is a \dmatroid.
\end{Proof}

Theorem~\ref{thm:vfclosed_del} suggests looking for an
excluded-minor characterization for the class of \vfsafe
matroids. By computer we found that the matroids $U_{2,6}$,
$U_{4,6}$, $P_6$, $F_7^-$, and $(F_7^-)^*$ (see again
\cite{Oxley/MatroidBook-2nd} for a description of these
matroids) are excluded minors for the family of \vfsafe
matroids. Moreover, using the database of D. Mayhew and G.F.
Royle \cite{DBLP:journals/jct/MayhewR08}, we confirmed that
these are the only excluded minors with $9$ or less elements.
We notice similarity with the excluded-minor characterization
of quaternary matroids (i.e., the matroids representable over
$GF(4)$) \cite{DBLP:journals/jct/GeelenGK00}, where it is shown
that a matroid $M$ is quaternary iff no minor of $M$ is
isomorphic to $U_{2,6}$, $U_{4,6}$, $P_6$, $F_7^-$,
$(F_7^-)^*$, $P_8$, or $P^=_8$. Hence, we conjecture the
following (which consequently has been verified for matroids
with $9$ or less elements).

\begin{Conjecture}
Every quaternary matroid is \vfsafe.
\end{Conjecture}

Let $\mathcal{N}$ be the family of matroids that have no minors
isomorphic to $U_{2,6}$, $U_{4,6}$, $P_6$, $F_7^-$, or
$(F_7^-)^*$. It is shown in
\cite[Corollary~1.2]{DBLP:journals/jct/GeelenOVW00} that
$\mathcal{N}$ can be constructed by taking direct sums and
2-sums of copies of $P_8^=$, minors of $S(5,6,12)$, and
quaternary matroids (see again \cite{Oxley/MatroidBook-2nd} for
a description of $S(5,6,12)$). In this light, we have verified
(by computer and using internal symmetries of the matroid) and
found that $S(5,6,12)$ is vf-safe as well. Hence we conjecture
that the above list of excluded minors for the class of \vfsafe
matroids is complete.

We finally note that not every quaternary \dmatroid is \vfsafe.
For example, the non-vf-safe \dmatroid $(V,2^V \setminus
\{\{u\}\})$ with $V = \{u,v,w\}$ (this \dmatroid differs from
the non-vf-safe \dmatroid of Section~\ref{sec:distance_DM} by a
pivot) is represented by the following skew-symmetric matrix
over $GF(4)$:
$$
\bordermatrix{& u & v & w\cr
u & 0 & a & b \cr
v & a & 1 & 0 \cr
w & b & 0 & 1
}
$$
where $a$ and $b$ are the two elements distinct from $0$ and $1$ in $GF(4)$.

\subsection*{Acknowledgements}
We thank Lorenzo Traldi for a stimulating correspondence. We
are much indebted to the anonymous reviewers for valuable comments
on earlier versions of this paper. We thank Gordon Royle for kindly
sending to us the database of matroids with nine elements from \cite{DBLP:journals/jct/MayhewR08}.
R.B.\ was supported by the Netherlands Organization for Scientific Research (NWO), project ``Annotated graph mining''.

\bibliography{nullity_loopc_dmatroid}

\end{document}